\documentclass[runningheads]{llncs}
\usepackage{algorithm}
\usepackage{amsmath}
\usepackage{amssymb}
\usepackage{algpseudocode}
\usepackage{graphicx}
\usepackage{hyperref,xcolor}
\usepackage{makeidx}
\definecolor{winered}{rgb}{0.5,0,0}
  
\hypersetup
{
    pdfauthor={Anurag Dutta},
    pdfsubject={Hyperlinks in LaTeX},
    pdftitle={myhyper.tex},
    pdfkeywords={LaTeX, PDF, hyperlinks}
    pdfborder={0 0 0},
colorlinks=true,
 linkcolor=winered,
urlcolor={winered},
filecolor={winered},
  linktoc=all,
}
\begin{document}
\title{A Comparative Investigation into the Operation of an Optimal Control Problem: The Maximal Stretch}
\titlerunning{The Maximal Stretch}
%
\author{Anurag Dutta\inst{1, *} \and
K. Lakshmanan\inst{2} \and
John Harshith\inst{3} \and
A. Ramamoorthy\inst{4}}
\authorrunning{A. Dutta et al.}
%
\institute{
Department of Computer Science, Government College of Engineering and Textile Technology, Serampore, Calcutta, India
\and 
Department of Mathematics, Kuwait American School of Education, Salmiya, Hawalli, Kuwait
\and  
Department of Computer Science, Vellore Institute of Technology, Vellore Campus, Vellore, Tamil Nadu, India 
\and
Department of Mathematics, Velammal Engineering College, Anna University, Chennai, Tamil Nadu, India
\\ \textbf{\\}
\email{anuragdutta.research@gmail.com}\\
}
\maketitle              

\begin{abstract}
Mathematical Selection is a method in which we select a particular choice from a set of such. It have always been an interesting field of study for mathematicians. Combinatorial optimisation is the practice of selecting the best constituent from a collection of prospective possibilities according to some particular characterization. In simple cases, an optimal process problem encompasses identifying components out of a finite arrangement and establishing the function's significance in possible to lessen or achieve maximum with a functional purpose. To extrapolate optimisation theory, it employs a wide range of mathematical concepts. Optimisation, when applied to a variety of different types of optimization algorithms, necessitates determining the best consequences of the specific predetermined characteristic in a particular circumstance. In this work, we will be working on one similar problem - The Maximal Stretch Problem with computational rigour. Beginning with the Problem Statement itself, we will be developing numerous step - by - step algorithms to solve the problem, and will finally pose a comparison between them on the basis of their Computational Complexity. The article entails around the Brute Force Solution, A Recursive Approach to deal with the problem, and finally a Dynamically Programmed Approach for the same. 
\keywords{Dynamic Programming \and Matrices \and Recursion \and Operation Research}
\end{abstract}
\section{Introduction}
Determining the contiguous subspace [1] with the biggest sum, $\sum_{i=\alpha}^{\beta}\mathcal{A}\left[i\right] \forall\mathcal{A}[i]\in \mathbb{R}$ together within specified linear array [2] $\mathcal{A}[1..n]$ of data, is defined as the "\textit{Maximum Summation SubArray Problem}" [3] in Computing Science. This problem is also referred to as the "\textit{Maximal Segmentation Sum Problem}". Ulf Grenander presented the maximal subarray problem [4] in 1977 as a streamlined framework for the expectation - maximization estimate of features in digitised images. In this article, we intend to provide an overview, and indeed outline a comprehensive evaluation, of the multiple alternatives to solving the problem - "The Maximal Stretch." The problem statement is as follows. \\\textbf{\\}
\textbf{Problem Statement: }\textit{Given a binary square matrix, $\mathcal{A}$ with binary entries, i.e, 
\begin{align*}
\mathcal{A}=\left[\begin{matrix}\alpha_{0,\ 0}&\cdots&\alpha_{0,\ n-1}\\\vdots&\ddots&\vdots\\\alpha_{n-1,\ 0}&\cdots&\alpha_{n-1,\ n-1}\\\end{matrix}\right]\ni\alpha_{i,\ j}=\left\{0,\ 1\right\}\ \ \forall\left(i,\ j\right)
\end{align*}
We will have to fetch a sub matrix
\begin{align*}
\mathcal{B}=\left[\begin{matrix}\beta_{0,\ 0}&\cdots&\beta_{0,\ m-1}\\\vdots&\ddots&\vdots\\\beta_{m-1,\ 0}&\cdots&\beta_{m-1,\ m-1}\\\end{matrix}\right]\ni\beta_{i,\ j}=1\ \ \forall\left(i,\ j\right)
\end{align*}
and find such a value of $m^2$ such that $\left\lfloor\frac{n}{m}\right\rfloor$ is minimized.}
\\\\A sparse matrix [5], also known as a sparse array, is a matrix where the majority of the entries are zero in quantitative simulation [6] and computer sciences. The total number of non-zero representatives is usually equivalent to the total number of the columns and rows, which is a common rule of thumb to determine whether such a structure is sparse. Nevertheless, there's currently no technical definition of the considered necessary percentage of zero-value modules. The matrix, on the other contrary, is considered to be intense if the significant proportion of its modules are non-zero.
\section{Naïve Approach}
The Naive Approach [7] is the one, which one's mind would spill out instantaneously after apprehending the Problem Statement. Infact, for the Maximal Stretch Problem, the Naïve Approach is mentioned in the Problem Statement itself, that is, look out for all sub matrices and find out the one having highest order with all affirmative entries. Generically, suppose, we have a matrix $\mathcal{M}$ of order $m \times n$, and let us declare a function $\varphi(m,\ n)$ which is ought to find out the number of sub - matrices [8] possible for $\mathcal{M}$. Let us try to build some recursive [9] stuff out of it, \\
\textit{Base Case}
\begin{equation}
\varphi(1,\ \lambda)=\varphi(\ \lambda,\ 1)=\sum_{i=1}^{\lambda}i
\end{equation}
\textit{Recursive Case\\}
Taking into consideration, the Matrix, $\mathcal{M}$ of order $m \times n$, there will be $\sum_{i=1}^{m}i$ sub - matrices of width 1. For each such matrix, there will be $n + 1$ sub - matrices. Hence, number of nascent sub - matrices will be $\left(n+1\right)\sum_{i=1}^{m}i$. So, the recursive form would be 
\begin{equation}
\varphi(m,\ n)=\varphi(m,\ n-1)+n\sum_{i=1}^{m}i
\end{equation}
The Closed form of this Recursion [10] comes out to be 
\begin{align*}
\varphi(m,\ n)=\left(\sum_{i=1}^{n}i\right)\times\left(\sum_{i=1}^{m}i\right)
\end{align*}
In our case, for the Maximal Stretch Problem, we have considered the Matrices to be Square in nature, so, $n = m$. Hence, the closed form of this Recursion [11] comes out to be
\begin{align*}
\varphi(n,\ n)=\sum_{i=1}^{n}i^3=\frac{n^2\left(n+1\right)^2}{4}
\end{align*}
Hence, the Computational Complexity for this Naïve Approach would be of Quartic [12] Order. The Algorithm corresponding to this Solution is mentioned in the Algorithm \ref{alg1}.
\begin{algorithm}
\caption{Naïve Solution to the Maximal Stretch Problem}\label{alg1}
\begin{algorithmic}[1]
\Require $\mathcal{M}_{n\times n}\ni\mathcal{M}_{\left(i,\ j\right)}=\left(0,\ 1\right)\ \forall\ \left(i,\ j\right)<\left(n,\ n\right)$
\Ensure $m^2$ $\exists \ {\mathcal{M}^\prime}_{m\times m}\ni{\mathcal{M}^\prime}_{\left(i,\ j\right)}=1\ \forall\ \left(i,\ j\right)<\left(m,m\right)$
      \Function{Naïve Solution}{$\mathcal{M}_{n\times n}$}
      \State $m \gets 0$ \Comment{$m$ stores the Maximal Length}
      \For{${\mathcal{M}^\prime}$ in $\mathcal{M}$}
                \If{${\mathcal{M}^\prime}_{\left(i,\ j\right)} = 0$}
                    \State \textbf{continue}
                \EndIf
                \State $m \gets max(m, $ ORDER$({\mathcal{M}^\prime}))$
      \EndFor
      \EndFunction
\end{algorithmic}
\end{algorithm}
\section{Recursive Approach}
Recursion is a mathematical problem-solving technique used in computing science in which the result is dependent on resolutions to local solutions of the exact same problem. To confront those very recurring issues, recursion employs practises that replicate themselves from within their own framework. Recursion is indeed a fundamental premise in computing science, and it can be used to alleviate a wide range of problems. Using a recursive approach, we will make an effort to achieve the consequence. The Methodology we will use is still quite simple.
\begin{enumerate}
\item Iterate over all elements of the matrix $\mathcal{M}$.
\item The recurrence subjecting to the problem can be declared as 
\begin{equation}
\psi\left(i,\ j\right)=1+min\left(\psi\left(i^+,\ j\right),\psi\left(i,\ j^+\right),\psi\left(i^+,\ j^+\right)\right)
\end{equation}
where, $i^+ = i + 1$ and $j^+ = j + 1$ with a base case of 
\begin{equation}
\psi\left(i,\ j\right)=0
\end{equation}
if ${\mathcal{M}}_{\left(i,\ j\right)}=0$
\item The maximum value [13] of $\psi\left(i,\ j\right)$ will give the value of $m$ for which $\left\lfloor\frac{n}{m}\right\rfloor$ is maximized.
\end{enumerate}
The Algorithm corresponding to this Solution is mentioned in the Algorithm \ref{alg2}.
\begin{algorithm}
\caption{Recursive Solution}\label{alg2}
\begin{algorithmic}[1]
\Require $\mathcal{M}_{n\times n}\ni\mathcal{M}_{\left(i,\ j\right)}=\left(0,\ 1\right)\ \forall\ \left(i,\ j\right)<\left(n,\ n\right)$
\Ensure $m^2$ $\exists \ {\mathcal{M}^\prime}_{m\times m}\ni{\mathcal{M}^\prime}_{\left(i,\ j\right)}=1\ \forall\ \left(i,\ j\right)<\left(m,m\right)$
      \Function{$\psi\left(i,\ j\right)$}{}
                \If{${\mathcal{M}}_{\left(i,\ j\right)} = 0$ or $\left(i,\ j\right)$ out of Range}
                    \State $\psi\left(i,\ j\right)=0$
                \Else
                    \State $\psi\left(i,\ j\right)=1+min\left(\psi\left(i+1,\ j\right),\psi\left(i,\ j+1\right),\psi\left(i+1,\ j+1\right)\right)$
                \EndIf
      \EndFunction
      \Function{Recursive Solution}{$\mathcal{M}_{n\times n}$}
      \State $m \gets 0$ \Comment{$m$ stores the Maximal Length}
      \For{$i=0$ to $n$}
      \For{$j=0$ to $n$}
      	\State $m \gets max(m, \psi\left(i,\ j\right))$
      \EndFor
      \EndFor
      \EndFunction
\end{algorithmic}
\end{algorithm}
The Computational Complexity of this Recursive Solution is of the order $n^2 \times \mathcal{T}_{\psi\left(i,\ j\right)}$. It is quite obvious that, $\mathcal{T}_{\psi\left(i,\ j\right)}=3^{n+n}$. since in each step, 3 computations are being performed. So, on a whole, the Computational Complexity of the problem is of the order $n^23^{2n}$.
\section{Dynamic Programming}
The significant proportion of languages used for computer programming endorse recursion by allowing a methodology to enact selves from within its own framework. A component that is called recurrently from the inside of itself may cause the request heap to develop to the dimensions of the system's total dimensions of all calls implicated. This means that recursion is typically less effective for challenges that can be solved quickly through iteration [14], and that using objective functions like dynamic programming [15] is essential for solving big problems. In this approach, we will make use of Dynamic Programming. The Algorithm, that will be followed follows as
\begin{enumerate}
\item We will create a 2 - Dimensional Bottom to Top look up table, 
$\mathcal{M}^\prime$ and initialize it with zeroes.
\item Iterate over the entries of the matrix $\mathcal{M}$
\begin{enumerate}
\item If $\mathcal{M}_{\left(i,\ j\right)} = 1$, ${\mathcal{M}^\prime}_{\left(i,\ j\right)}=min({\mathcal{M}^\prime}_{\left(i-1,\ j\right)},{\mathcal{M}^\prime}_{\left(i,\ j-1\right)},\ {\mathcal{M}^\prime}_{\left(i-1,\ j-1\right)})+1\ $.
\item Else ${\mathcal{M}^\prime}_{\left(i,\ j\right)}=0$.
\end{enumerate}
\item Return the maximum value of ${\mathcal{M}^\prime}_{\left(i,\ j\right)}$.
\end{enumerate}
The Algorithm corresponding to this Solution is mentioned in the Algorithm \ref{alg3}.
\begin{algorithm}
\caption{DP Solution}\label{alg3}
\begin{algorithmic}[1]
\Require $\mathcal{M}_{n\times n}\ni\mathcal{M}_{\left(i,\ j\right)}=\left(0,\ 1\right)\ \forall\ \left(i,\ j\right)<\left(n,\ n\right)$
\Ensure $m^2$ $\exists \ {\mathcal{M}^\prime}_{m\times m}\ni{\mathcal{M}^\prime}_{\left(i,\ j\right)}=1\ \forall\ \left(i,\ j\right)<\left(m,m\right)$
      \Function{DP Solution}{$\mathcal{M}_{n\times n}$}
      \State ${\mathcal{M}^\prime}_{(n+1)\times (n+1)} \gets 0$
      \State $m \gets 0$ \Comment{$m$ stores the Maximal Length}
      \For{$i=0$ to $n$}
      \For{$j=0$ to $n$}
                \If{${\mathcal{M}}_{\left(i,\ j\right)} = 0$}
                    \State ${\mathcal{M}^\prime}_{\left(i,\ j\right)}=0$
                \Else
                    \State ${\mathcal{M}^\prime}_{\left(i,\ j\right)}=min({\mathcal{M}^\prime}_{\left(i-1,\ j\right)},{\mathcal{M}^\prime}_{\left(i,\ j-1\right)},\ {\mathcal{M}^\prime}_{\left(i-1,\ j-1\right)})+1\ $
                \EndIf
                \State $m \gets max(m, {\mathcal{M}^\prime}_{\left(i,\ j\right)})$
      \EndFor
      \EndFor
      \EndFunction
\end{algorithmic}
\end{algorithm}
The Computational Complexity for this DP Solution is of the Quadratic Order.
\section{Conclusion}
To Conclude our work, The Maximal Stretch Problem can be solved in 3 ways 
\begin{enumerate}
\item By looking up all possible Sub matrices,
\begin{align*}
\mathcal{B}=\left[\begin{matrix}\beta_{0,\ 0}&\cdots&\beta_{0,\ m-1}\\\vdots&\ddots&\vdots\\\beta_{m-1,\ 0}&\cdots&\beta_{m-1,\ m-1}\\\end{matrix}\right]\ni\beta_{i,\ j}=1\ \ \forall\left(i,\ j\right)
\end{align*} 
that are possible for a given matrix, 
\begin{align*}
\mathcal{A}=\left[\begin{matrix}\alpha_{0,\ 0}&\cdots&\alpha_{0,\ n-1}\\\vdots&\ddots&\vdots\\\alpha_{n-1,\ 0}&\cdots&\alpha_{n-1,\ n-1}\\\end{matrix}\right]\ni\alpha_{i,\ j}=\left\{0,\ 1\right\}\ \ \forall\left(i,\ j\right)
\end{align*}
Now, this approach will compute the result in Quartic Order, $n^4$.
\item By Recurring though all possible indices of the matrix, for the relation
\begin{align*}
\psi\left(i,\ j\right)=\left\{\begin{matrix}1+min\left(\psi\left(i^+,\ j\right),\psi\left(i,\ j^+\right),\psi\left(i^+,\ j^+\right)\right)\\0\\\end{matrix}\right.
\end{align*}
This will compute the result in a Hyper tonic order of $n^23^{2n}$.
\item By Preparing a Look Up table to implement the recurrence mentioned above. This will compute the result in Quadratic Time, $n^2$.
\end{enumerate}
\begin{figure}[htbp]
\centerline{\includegraphics[width = \linewidth]{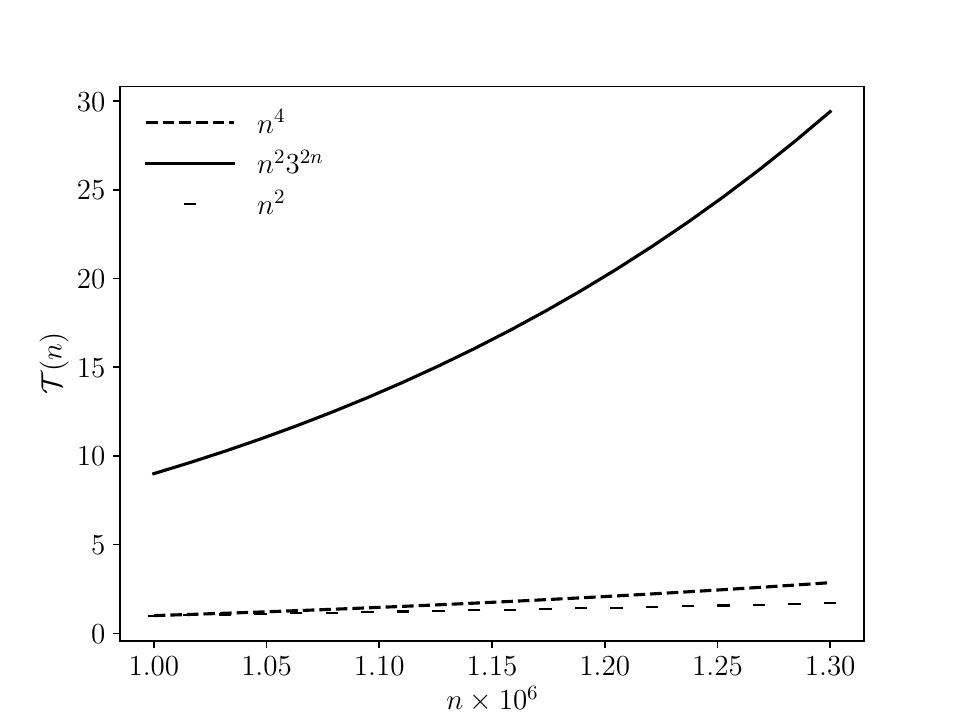}}
\label{fig1}
\caption{Comparative Plot between the 3 Techniques for solving the Maximal Stretch Problem. The computational time is taken on the Y - axes, and the Cardinality of the Dataset subjected to the problem is taken on the X - axes.}
\end{figure}
It is quite observable from the Comparative Plot, that the Computational Algorithm of the problem subjective to Dynamic Programming turns out to be the best. 

\end{document}